\documentclass[12pt]{amsart}

\usepackage{amsmath}
\theoremstyle{plain}
\newtheorem{theorem}{Theorem}
\newtheorem{corollary}[theorem]{Corollary}
\newtheorem{lemma}[theorem]{Lemma}
\newtheorem{proposition}[theorem]{Proposition}

\theoremstyle{definition}
\newtheorem{definition}[theorem]{Definition}

\theoremstyle{remark}

\begin{document}
\def\complex{\bf C}
\def\R{{\rm I\kern-2ptR}}
\def\N{{\rm I\kern-2ptN}}
\def\complex{{\bf C}}
\def\C{{\bf C}}
\def\CC{{\footb C}}
\def\blambda{\mbox{\boldmath$\lambda$}}
\def\bomega{\mbox{\boldmath$\omega$}}
\def\emptyset{\mathop{\raisebox{.1ex}{$\not
\mathrel{\raisebox{.1ex}{$\scriptstyle\bigcirc$}}$}}}

\def\qed{\hskip .6em \raise1.8pt\hbox{\vrule height4pt
width6pt depth2pt}}
\def\qedd{\hskip .4em \raise1.8pt\hbox{\vrule height3pt
width5pt depth1.8pt}}
\def\sq{\hskip .6em \raise1.8pt\hbox{\vrule
height4pt width6pt depth2pt}}

\def\lbr{\bigl[\hskip-.06cm[}
\def\rbr{\bigr]\hskip-.06cm]}

\def\od#1{#1^\sim}
\def\pref{\succeq}
\def\spref{\succ}
\def\indiff{\sim}
\def\notpref{\not\pref}
\def\s#1{{#1^*}}

%\font\Bbb=msbm10
%\def\R{\hbox{\Bbb R}}
%\def\N{\hbox{\Bbb N}}

\font\titlos=cmb10 at 18 pt
%\font\titlos=cmb10 at 20 pt
\font\head=cmbx8 at 33 pt
\font\titlosa=cmbx9 at 18 pt

%\hfill DRAFT, \ \ \today

\medskip

\centerline{\titlos A Characterization of Compact-friendly}
\vskip .12cm 
\centerline{\titlos Multiplication Operators}

\vskip 1.3cm

\centerline{\bf Y.~A. Abramovich$^a$, C.~D. Aliprantis$^b$, 
O.~Burkinshaw$^a$ and A.~W. Wickstead$^c$}

\vskip .4cm
\noindent{$^a$\footnotesize Department of Mathematical Sciences,
IUPUI, Indianapolis, IN 46202-3216, USA}

\noindent{$^b$\footnotesize Departments of Economics and  Mathematics,
Purdue University, West Lafayette, IN 47907, USA}

\noindent{$^c$\footnotesize Department of Pure  Mathematics,
  The Queen's University of Belfast,  Belfast BT7 1NN, UK}

\vskip 1.5cm

\centerline{
\vbox{\hsize 5 truein \baselineskip=12pt
\footnotesize 
\noindent ABSTRACT. Answering in  the affirmative a question posed 
in~\cite{AAB3}, we prove that a positive multiplication 
operator on any $L_p$-space (resp. on a $C(\Omega)$-space) 
is compact-friendly if and only if the  multiplier is constant 
on a set of positive measure (resp. on a non-empty
open set).\hfill\break\indent  
    In the process of establishing this result, we also prove that
any  multiplication operator has a family of hyperinvariant 
bands---a fact  that does not seem to have appeared in the 
literature before. This provides useful information about 
the commutant of a multiplication operator.
}}

\vskip 1.1cm

\section{Preliminaries}\label{sec1}    %\marginnote{sec1}

This work will employ techniques and terminology from Banach lattice theory.  For
terminology which is not explained below, we refer the
reader to~\cite{AB}.

In this work the word ``operator'' will be synonymous 
with ``linear operator.'' An operator $T\colon X\to Y$ between two
Banach lattices is  {\it positive\/} if $x\ge0$ in $X$ implies $Tx\ge0$ in $Y$. 

A positive operator $S\colon X\to X$ on a Banach lattice $X$ is said to
{\bf dominate} another operator $T\colon X\to X$ (in symbols, $S\succ T$) if 
\[
|Tx|\le S|x|
\]
for each $x\in X$.  If $S$ dominates $T$, we shall also say that $T$ {\it is
dominated\/} by $S$. Every operator dominated by a positive operator is
automatically continuous. 

We recall next the notion of a compact-friendly operator that was introduced 
in~\cite{AAB1} and that will play an important role in this work.

\begin{definition}\label{def:Com-fr}      %\marginnote{def:Com-fr}
A positive operator $B\colon X\to X$ on a Banach lattice        
is said to be {\bf compact-friendly} if there exist three non-zero
operators $R, K, A\colon X\to X$ with $R$ and $ K$  positive and $K$
compact  satisfying
$$
RB=BR,\quad R\succ A\quad\hbox{and}\quad K\succ A\,.
$$
\end{definition}

Regarding the invariant subspace problem for operators on Banach lattices the
compact-friendly operators seem to be the analogues of Lomonosov operators. 
Recall that an operator $T\colon X\to X$ on a Banach space is  a    
{\bf Lomonosov operator} if there exist non-zero operators $S,K\colon X\to X$ 
such that {\it $S$ is not a multiple of the identity, $K$ is compact,
$ST=TS$, and $SK=KS$.\/}

The invariant subspace theorems for positive operators obtained in~\cite{AAB1}
(see also~\cite{AAB2}) can be viewed as the Banach lattice analogues of the
following famous invariant subspace theorem of V.~I. Lomonosov.

\begin{theorem}[Lomonosov~\cite{LO}]\label{Lomon}      %\marginnote{Lomon}
Every Lomonosov operator $T$  has a non-trivial clo\-sed invariant subspace.
Moreover, if $T$  itself commutes with a non-zero compact operator,
then there exists  a  non-trivial closed hyperinvariant subspace.
\end{theorem}

Besides compact-friendly operators, we shall work here also with
multiplication operators on spaces of continuous and measurable functions.      
If $\Omega$ is a compact Hausdorff space and $\phi\in C(\Omega)$, then 
a {\bf multiplication operator} $M_\varphi$ on  $C(\Omega)$ is  defined by
$M_\varphi f=\phi f$ for each $f\in C(\Omega)$. The function $\phi$ is called  
the {\it multiplier\/}.

    Similarly, if $X$ is  a  Banach function space on a measure space
$(\Omega,\Sigma,\mu)$ and  $\phi$ is a measurable function,
then a multiplication operator $M_\varphi$ on $X$ is defined by
$M_\varphi f=\phi f$ for each $f\in X$. 
Observe that  a multiplication operator $M_\phi$ maps $X$ into itself           
if and only if the multiplier $\phi$  is an (essentially) bounded function.
So, for the rest of this paper, whenever we deal with a multiplication    
operator $M_\phi$ on a Banach function space we assume that the multiplier 
$\phi\in L_\infty(\mu)$.

It should be noticed that a multiplication operator is positive if and only     
if its multiplier is a non-negative function. 

Obviously  each multiplication operator $M_\varphi$ has non-trivial 
invariant subspaces and, as was observed in~\cite{AAB3}, 
each multiplication operator  is a Lomonosov operator. 
Moreover, as we
will prove in the next section (see Theorem~\ref{strips}  and   
Corollary~\ref{c:strips}) each multiplication operator $M_\varphi$ has
hyperinvariant subspaces of a very simple geometrical form, namely,            the
disjoint bands.

Our next definition describes the kind of multipliers that will
be important in our work.

\begin{definition}\label{d:flat}        %\marginnote{d:flat}
A continuous function $\phi\colon\Omega\to\R$ on a topological space has a
{\bf flat} if there exists a non-empty open set $V$ such that $\phi$ is
constant on $V$.

Similarly, a measurable function $\phi\colon \Omega\to\R$ on a
measure space $(\Omega,\Sigma,\mu)$ is said to have a {\bf flat} if 
$\phi$ is constant on some $A\in\Sigma$ with $\mu(A)>0$.
\end{definition}

It was shown in~\cite{AAB3} that a positive multiplication operator commutes  
with a non-zero finite rank operator if and only if the multiplier has a flat. 
It was then asked whether the flatness condition  characterizes also the 
compact-friendly multiplication operators. The objective of this work is   to
answer this question affirmatively.  Namely, the main result of this paper can  
be stated as follows.

\begin{theorem}\label{th:MAIN}    %\marginnote{th:MAIN}
A positive multiplication operator $M_\varphi$ on a $C(\Omega)$-space or on a
$L_p(\Omega,\Sigma,\mu)$-space $(1\le p\le \infty)$ is compact-friendly if and 
only if the multiplier $\varphi $ has a flat.  
\end{theorem}

\section{The commutant of a multiplication operator}\label{sec2}

In this section $X$ will denote a Banach function space on a fixed  measure
space $(\Omega,\Sigma,\mu)$. Let $M_\phi\colon X\to X$ be the 
multiplication operator with a multiplier $\varphi \in L_\infty(\mu)$.

Not much is known about the commutant of $M_\varphi$. The following
discussion will provide some important insights
into the structure of the commutant. We precede this discussion by  
fixing some notation. If $f\colon\Omega\to\R$ is a function, then its   
support,   ${\rm Supp}(f)$, is defined by
\[
{\rm Supp}(f)=\{\omega\in\Omega\colon\ f(\omega)\neq0\}\,.
\]
If $A,B\in\Sigma$, then relations $A\subseteq B$ a.e. and $A=B$ a.e.   
are understood as usual $\mu$-a.e. For example, $A\subseteq B$ a.e. means
that  $\mu(\{\omega\in A\colon\ \omega\notin B\})=0$.

\begin{definition}\label{d:suppt}     %\marginnote{d:suppt} 

Let $T\colon X\to X$ be a continuous operator and let $E\subseteq \Omega $     
be a measurable subset of positive measure. We shall say  that $T$            
{\bf leaves $E$ invariant}, if
\[
x\in X \ \hbox{ and }\ {\rm Supp}(x)\subseteq E\ \hbox{ implies }\
{\rm Supp}(Tx)\subseteq E \ a.e.
\] 
\end{definition}
This definition is, of course, motivated by a simple observation
that an operator $T$ leaves a
(measurable) set $E$ invariant if and only if $T$ leaves invariant the band
$$
B_E=\{f\in X\colon\ f=0\,\hbox{ on }\,\Omega\setminus E\}
$$ 
generated by $E$ in $X$. It is obvious that {\it if an operator           
$T\colon X\to X$ leaves invariant the sets $E$ and $F$, then it also        
leaves invariant the sets $E\cap F$ and $E\cup F$.\/}

Now let us introduce some more notation. 
For each  $\alpha\in \R$, let
$$
E_\alpha = \bigl\{\omega\in\Omega\colon\ \phi(\omega) \ge
\alpha\bigr\}\quad\hbox{and}\quad
E^\alpha =\bigl\{\omega\in\Omega\colon\ \phi(\omega) \le
\alpha\bigr\}\,.
$$
If we need to emphasize that the level set $E_\alpha$ is produced by the
function $\phi$, then we shall write $E_\alpha(\phi)$ instead of $E_\alpha$.
For $\alpha \le\beta$,  we also write
$$
E_\alpha^\beta =E_\alpha\cap E^\beta=\bigl\{\omega\in\Omega\colon\ \alpha \le
\phi(\omega) \le \beta\bigr\}\,. 
$$

And now we come to a simple but important result asserting that all the bands
in $X$ generated by   the  level sets  introduced above are left invariant 
by each operator commuting with $M_\phi$.

\begin{theorem}\label{strips}      %\marginnote{strips}

Every operator in the commutant of $M_\phi$ leaves invariant all the 
sets $E_\alpha$, $E^\alpha$ and $E_\alpha^\beta$.    
\end{theorem}

\begin{proof}
Let $R\colon X\to X$ be a bounded operator commuting with $M_\varphi$. 
We begin by considering the sets $E^\alpha$.
Assume  that $\varphi \ge 0$. First we will verify that $R$
leaves  invariant the set $E^\alpha$ with $\alpha=1$,  i.e., the set
$$ 
E^1=E^1(\varphi)=\{\omega\in\Omega\colon\ \varphi(\omega)\le1\}.
$$  
To do this, assume by way of contradiction that $R$ does not leave $E^1$
invariant. This means that there exists some function $x\in X$ with             
${\rm Supp}(x)\subseteq E^1$ and such that the measurable set
$A=\{\omega\in\Omega\colon\  Rx(\omega)\neq0\ \& \ \varphi(\omega)>1\}$            
has positive measure. Pick some $\gamma >1$ such that  
$B=\{\omega\in\Omega\colon\  Rx(\omega)\neq0\ \& \
\varphi(\omega)>\gamma\}$ has positive measure.  

The commutativity property $RM_\varphi=M_\varphi R$ easily implies
$$ 
R(\varphi^nx)=\varphi^nRx\eqno(\star)
$$ 
for each $n$. Let $\|\cdot\|$ denote the  norm on $X$. We shall reach a
contradiction by computing the norm of the function in
$(\star)$ in two different ways. On  one hand, the hypothesis 
${\rm Supp}(x)\subseteq E^1$ and the fact that $0\le \varphi(\omega) \le 1$ on  
$E^1$ imply that $|\varphi^nx|\le|x|$, and so
$$
\|R(\varphi^nx)\|\le\|R\|\cdot\|\varphi^nx\|\le\|R\|\cdot\|x\|<\infty\,.
$$ 
On the other hand, for the element $y=|(Rx)\chi_B| \in X$ we have
$$ 
0<\gamma^n y\le|\varphi^n (Rx) \chi_B|\le|\varphi^nRx|\,,
$$  
whence
$$ 0< \gamma^n\|y\|\le \|\varphi^nRx\|=
\|R(\varphi^nx)\|\le\|R\|\cdot\|x\|<\infty
$$ 
for each $n$, contradicting the fact that $\gamma>1$. Hence, $R$ leaves
$E^1(\varphi)$ invariant. 

Let us verify now that $R$ leaves invariant each $E^\alpha$ with
$\alpha > 0$. Consider $\psi=\alpha^{-1} \varphi$. Obviously 
the multiplication operator $M_\psi$ also commutes with $R$ and $E^1(\psi)=
E^\alpha(\varphi)$. By the previous part $R$ leaves $E^\alpha$ invariant.

Since $E^0 =\cap_{\alpha> 0}E^\alpha$ we see that $R$ leaves
$E^0$ invariant as well.  Since $\varphi \ge 0$ the set
$E^\alpha =\emptyset$ whenever $\alpha < 0$.
Thus, for $\varphi \ge 0$ we have proved that
$R$ leaves any set $E^\alpha$ invariant. The assumption
made at the beginning of the proof that the multiplier $\varphi$ 
is nonnegative can be easily  disposed of. Indeed, pick any  
$t>0$ such that the function          
$\psi=\varphi +t {\bf 1}$ is positive.
Obviously $M_\psi$ commutes with $R$ (since
$M_\varphi$  does) and $E^\alpha(\varphi)= E^{\alpha+t}(\psi)$. 
By the preceding part $R$ leaves $E^{\alpha+t}(\psi)$, that is 
$E^\alpha(\varphi)$, invariant.

Finally notice  that $E_\alpha(\varphi)=E^{-\alpha}(-\varphi)$.
This shows that the case of the sets $E_\alpha$ follows immediately from
the case of the sets $E^{\alpha}$ consided above. 
\end{proof}

\begin{corollary}\label{c:strips}      %\marginnote{c:strips}
If $\phi\in L_\infty(\mu)$ is a non-constant function, then 
the  multiplication operator $M_\phi$ has a non-trivial hyperinvariant       
band. If the $($essential$)$ range of the multiplier $\varphi $ is 
an infinite set, then $M_\varphi$ has infinite many disjoint 
hyperinvariant bands.
\end{corollary}

Consider also    the
following three additional types of the level sets associated with
the multiplier $\varphi$: 
$$
\{\omega\in\Omega\colon\ \alpha \le \varphi(\omega) < \beta\}, \
\{\omega\in\Omega\colon\ \alpha < \varphi(\omega) < \beta\} \ \hbox{ and }\
\{\omega\in\Omega\colon\ \alpha < \varphi(\omega) \le \beta\}. 
$$

It is easy to see that if  $R$ is order continuous (and commutes with $M_\varphi$) 
then $R$ leaves also each of these sets invariant. In particular this is so
if the norm on $X$ is order continuous. However, quite surprisingly, it may    
happen that without this extra assumption the operator $R$  may fail to leave   
these latter sets invariant.

Even when $X$ has order continuous norm (and so  $R$ leaves   
invariant  so many mutually disjoint bands) it is not true in general that        
$R$ leaves  invariant any band.  
Furthermore, as we shall see in the next example 
$R$ may even fail to be a disjointness  preserving operator. (Recall that an 
operator $R$ on a vector lattice is said to  preserve disjointness if
$R$ carries disjoint vectors to disjoint vector.)

\begin{itemize}
\item A positive  operator  $R\colon L_\infty \to L_\infty $ commuting with
$M_\varphi$  need not be disjointness preserving  even if $\varphi$ has no flat. 
\end{itemize}

\noindent  To see this take $\mu$ to be the usual 2-dimensional Lebesgue 
measure on $[0,1]\times[0,1]$ and $\varphi(x,y)=y$. Let 
$Rf(x,y)=\int_0^1 f(t,y)\;dt$, then it is easy to see that $R$ 
commutes with $M_\varphi$, the multiplier $\varphi$ has no flat but $R$ is not
disjointness preserving. [If $\varphi$ has a flat, then the existence
of  $R$ as required is obvious].

\section{Multiplication operators on $C(\Omega)$-spaces}\label{sec4}
%\marginnote{sec4}

We start with a useful  general criterion 
for distinguishing between compact-friendly and non-compact-friendly         
operators on a Banach lattice with order continuous norm.

\begin{proposition}\label{t:general} %\marginnote {t:general}

Let $A\colon Y\to Y$ be an operator on a Banach lattice
dominated by a positive  compact operator. Then for any norm bounded sequence
$\{e_n\}$ the following two statements are true.
\begin{enumerate}
\item
The sequence $\{Ae_n\}$ has an order bounded subsequence.
\item
If $Y$ has order continuous norm and $\{Ae_n\}$ is disjoint, then
$\|Ae_n\| \to 0$. 
\end{enumerate}
\end{proposition}

\begin{proof}
(1) Let $K\colon Y\to Y$ be a compact positive operator dominating $A$, i.e.,
$|Ax|\le K|x|$ holds for each $x\in Y$. Since $K$ is a compact operator and
$\{e_n\}$ is a norm bounded sequence, we can extract from $\{K(|e_n|)\}$
a convergent subsequence. Without loss of generality we can assume that
the sequence $\{K(|e_n|)\}$ itself converges in $Y$, that is, there exists
$y\in Y$ such that $K|e_n|\to y$. By passing
to another subsequence if necessary, we can also assume without loss of
generality that $\|K|e_n|-y\|<2^{-n}$ holds for each $n$. Letting  
$e=\sum_{n=1}^\infty|K|e_n|-y|$ we see that $e\in Y^+$ and clearly  
$|K|e_n|-y|\le e$, whence $K|e_n|\le e+|y|$ for each $n$.              
It remains to note that
\[
|Ae_n|\le K|e_n|\le e+|y|
\]
for each $n$.
\vskip .12cm
(2) Assume that $\{Ae_n\}$ is a disjoint sequence and let $\{f_n\}$ be a
subsequence of $\{e_n\}$. By part (1), there exists a subsequence $\{g_n\}$
of $\{f_n\}$ (and hence of $\{e_n\}$) such that the pairwise disjoint     
sequence $\{Ag_n\}$ is order bounded.  Since $Y$ has order continuous norm,     
it follows that $Ag_n\to 0$ in $Y$; see~\cite[Theorem~{12.13}, p.~{183}]{AB}.
Thus, we have shown that every subsequence of $\{Ae_n\}$ has a 
subsequence convergent to zero, and consequently $Ae_n\to0$ in $Y$.
\end{proof}

The next theorem is a characterization of the compact-friendly 
multiplication operators on $C(\Omega)$-spaces.

\begin{theorem} \label{t:tony}  %\marginnote{t:tony}
A positive multiplication operator $M_\varphi$ on a $C(\Omega)$-space is
compact-friendly if  and only the multiplier $\varphi$ has a flat.
\end{theorem}
\begin{proof}
Let $0\le\phi\in C(\Omega)$. If  $\phi$ is constant on a non-empty open subset
of $\Omega$, then $M_\phi$ commutes with a non-zero positive rank-one
operator (see~\cite[Theorem~{2.6}]{AAB3}), and so $M_\varphi$ is compact-friendly.

For the converse, assume that $M_\phi$ is compact-friendly, and
consequently  there exist non-zero bounded operators $R,K,A\colon C(\Omega)\to
C(\Omega)$ with $R, K$ positive, $K$ compact and such that
$$
M_\phi R=RM_\phi,\ \ R\succ A\ \ \hbox{and}\ \ K\succ A\,.
$$
Taking adjoints, we see that
$$
M_\phi^* R^*=R^*M_\phi^*,\ \ R^*\succ A^*\ \ \hbox{and}\ \ 
K^*\succ A^*\,.
$$
The following three properties follow in a rather straightforward way.

1) For each $\omega\in \Omega$ the support of the measure 
$R^*\delta_\omega$ is contained in the set
$W_\omega=\phi^{-1}(\phi(\omega))$, 
where $\delta_\omega$ denotes the unit mass at $\omega$.
This claim is immediate from consideration of the identity
$$
M_\phi^* R^*\delta_\omega=R^*
M_\phi^*\delta_\omega=\phi(\omega)R^*\delta_\omega\,.
$$

2) Since $R\succ A$, it follows immediately from $1)$ that   
for each $\omega\in \Omega$ the measure $A^*\delta_\omega$ is also     
supported  by $W_\omega$.

3)  Pick $h\in C(\Omega)$ with $\|h\|=1$ and
$Ah\ne 0$. Next, choose a non-empty open set $U$ on which $|Ah(\omega)|\ge
\epsilon>0$ for some $\epsilon>0$. Then
for each $\omega\in U$ we have
$\|A^*\delta_\omega\|\ge \epsilon$.
Indeed, to see this, notice that
$$
\|A^*\delta_\omega\|\ge | \langle A^*\delta_\omega,h\rangle|=
|\langle \delta_\omega,Ah\rangle| =|Ah(\omega)|\ge \epsilon.$$

To complete the proof, assume by way of contradiction that
the set $W_\omega$ has an empty interior for each $\omega\in \Omega$.
Then the non-empty open set $U$, chosen in (3) must
meet infinitely  many sets $W_\omega$.  Pick a sequence $\{\omega_n\}$ in
$U$ with $\phi(\omega_m)\ne \phi(\omega_n)$ if $m\ne n$, 
and let $e_n=|A^*\delta_{\omega_n}|$ for each $n$. Then
$\|e_n\|\ge\epsilon$ for each $n$.  Furthermore,
since each $e_n$ is supported by the set $W_{\omega_n}$
and the sequence $\{W_{\omega_n}\}$ is pairwise disjoint, the sequence
$\{A^*\delta_{\omega_n}\}$ is also disjoint. However, by
Proposition~\ref{t:general} (which is applicable since the norm in $C(\Omega)^*$ 
is order continuous) we should have $\|A^*\delta_{\omega_n}\|\to 0$, a
contradiction. This completes the proof of the theorem. 
\end{proof}

Since each $L_\infty(\mu)$ space can be represented as $C(\Omega)$
space on its Stone space, the previous theorem implies immediately 
the following result.

\begin{theorem}\label{t:AM.const} %\marginnote {t:AM.const}
A multiplication operator $M_\phi$ on $L_\infty$, where             
$\phi\in L_\infty(\mu)$, is compact-friendly if and only if its
multiplier $\phi$ has a flat.
\end{theorem}

\section{Compact-friendly multiplication operators on $L_p$-spaces}\label{sec3}
%\marginnote{sec3}

For the rest of  our discussion, $(\Omega,\Sigma,\mu)$ 
will denote a fixed measure space, and $\|\cdot\|$ will denote the        
standard norm on $L_p(\mu)$. The main result in this section 
is the following  $L_p$-version of Theorems~\ref{t:tony} 
and~\ref{t:AM.const}.

\begin{theorem}\label{t:oc.const} %\marginnote {t:oc.const}
A multiplication operator $M_\phi$ on an arbitrary  $L_p(\mu)$-space, 
where   $0\le \phi\in L_\infty(\mu)$ and $1\le p<\infty$,
is compact-friendly if  and only if $\varphi$ has a flat.\footnote{We do
not know if this  theorem is  true for arbitrary Banach function
spaces.}   
\end{theorem}
\begin{proof}
It was shown in~\cite{AAB3} that if  $\phi$ has a flat, then $M_\phi$ 
commutes with a positive rank-one operator---and hence $M_\phi$ is 
compact-friendly. 

In the converse direction, assume that $M_\phi$ is compact-friendly and      
that, contrary to our claim, $\phi$ is not constant on any set
of positive measure. Pick three non-zero bounded operators 
$R,A,K\colon L_p(\mu)\to L_p(\mu)$ such that $R$ and $ K$ are
positive, $K$ is compact  and 
$$
RM_\phi=M_\phi R,\quad R\succ A\quad\hbox{and}\quad K\succ A\,.
$$
To obtain a contradiction, it will suffice (in view of
Proposition~\ref{t:general})   to construct a sequence  $\{e_n\}$ in $L_p(\mu)$
satisfying the following     properties:
\begin{itemize}
\item[(i)]
$\|e_n\| =1$ for each $n$,
\item[(ii)]
$\{Ae_n\}$ is a disjoint sequence, and
\item[(iii)]
$\|Ae_n\| \ge \delta$ for each $n$ and for some $\delta >0$.
\end{itemize}

The construction of such a sequence is quite involved and  will be presented 
in a series of lemmas below.
\end{proof}

The rest of this section will be devoted to construction of a sequence  
$\{e_n\}$ that satisfies the properties (i), (ii) and (iii) stated at the      
end of the proof of Theorem~\ref{t:oc.const}. We begin with some preliminary
comments. 
\begin{enumerate}
\item 
The assumption that $\phi$ does not have a flat means that for each   
$\gamma\ge0$ the set $E_\gamma^\gamma=
\{\omega\in\Omega\colon\ \phi(\omega)=\gamma\} =\phi^{-1}(\{\gamma\})
$
has measure zero. In particular, this implies that for any 
$\gamma \in (\alpha,\beta)$ the level sets $E_\alpha^\gamma$ and $
E_\gamma^\beta$    are essentially disjoint (in the sense that $E_\alpha^\gamma
\cap E_\gamma^\beta=E_\gamma^\gamma$ is  a set of measure zero).

\item
By Theorem~\ref{strips} the operator $R$ leaves  all the level sets of       
$\phi$ invariant, and so does the operator
$A$ since it is  dominated by $R$.

\item
Since $A\neq 0$ there exists some $x\in L_p(\mu)$ with $y=Ax\neq 0$. The
functions $x$ and $y$ will be fixed throughout the discussion in this section.
If we let $\alpha_0=0$ and $\beta_0=\|\phi\|_\infty$, then obviously 
$E_{\alpha_0}^{\beta_0}=\Omega$ and so 
\[
{\rm Supp}(x)\subseteq E_{\alpha_0}^{\beta_0}\,.
\] 
\end{enumerate}

\begin{lemma}\label{l:norm}  %\marginnote{l:norm}

There exists some $\gamma_0 \in (\alpha_0, \beta_0)$ such that 
$$
\bigl\|y\chi_{E_{\alpha_0}^{\gamma_0}}\bigr\| = 
\bigl\|y\chi_{E_{\gamma_0}^{\beta_0}}\bigr\| = c \|y\|\,,              
$$
where $c=1/{\root p \of 2}$.
\end{lemma}

\begin{proof}
Consider the function $N\colon[\alpha_0,\beta_0] \to \R$ defined by
$$
N(\gamma) = \|y\chi_{E_{\alpha_0}^\gamma}\|\,.
$$
Clearly, $N(\alpha_0)=0, N(\beta_0) =\|y\|$, and the function $N$ is   
continuous  by virtue of the ``no flats'' assumption about $\varphi$.   
Therefore, there exists  some 
$\gamma_0 \in (\alpha_0, \beta_0)$ such that $N(\gamma_0)=c \|y\|$.

Since $y\chi_{E_{\alpha_0}^{\gamma_0}}+ y\chi_{E_{\gamma_0}^{\beta_0}} =y$, 
and since the sets $E_{\alpha_0}^{\gamma_0}, E_{\gamma_0}^{\beta_0}$ are      
essentially disjoint, the $p$-additivity of the norm in $L_p(\mu)$ implies
that 
$$
\|y\chi_{E_{\alpha_0}^{\gamma_0}}\|^p +\|y\chi_{E_{\gamma_0}^{\beta_0}}\|^p 
=\|y\|^p.
$$
Consequently,
$$
\|y\chi_{E_{\gamma_0}^{\beta_0}}\|^p =\|y\|^p
-\|y\chi_{E_{\alpha_0}^{\gamma_0}}\|^p=\|y\|^p -c^p\|y\| =     
\frac{1}{2}\|y\|^p= c^p\|y\|^p,
$$
that is, $\|y\chi_{E_{\gamma_0}^{\beta_0}}\| =c\|y\|$, as required.
\end{proof}

Using the sets $E_{\alpha_0}^{\gamma_0}$ and $E_{\gamma_0}^{\beta_0}$
we can  represent $x$ as 
$$
x=x\chi_{E_{\alpha_0}^{\gamma_0}}\oplus
x\chi_{E_{\gamma_0}^{\beta_0}}\,, 
$$
and denote by $a_1$ the summand with smaller (or equal)
norm. The other summand will be denoted by $b_1$. So, if
$\|x\chi_{E_{\alpha_0}^{\gamma_0}}\|\le \|x\chi_{E_{\gamma_0}^{\beta_0}}\|$, 
then we let $a_1= x\chi_{ E_{\alpha_0}^{\gamma_0}}$ and 
$b_1= x\chi_{E_{\gamma_0}^{\beta_0}}$, and thus
$$
x= a_1\oplus b_1\,.
$$
Having chosen $a_1$ and $b_1$, we let
$$
                 u_1=y\chi_{E_{\alpha_0}^{\gamma_0}}\quad\hbox{and}\quad
                 v_1=y\chi_{E_{\gamma_0}^{\beta_0}}
$$
and also  $\alpha_1=\alpha_0$ and $\beta_1=\gamma_0$. 
(However, if 
$\|x\chi_{E_{\gamma_0}^{\beta_0}}\| <\|x\chi_{E_{\alpha_0}^{\gamma_0}}\|$,
then $a_1= x\chi_{E_{\gamma_0}^{\beta_0}}$, 
$b_1= x\chi_{E_{\alpha_0}^{\gamma_0}}$,
and we let $u_1=y\chi_{E_{\gamma_0}^{\beta_0}}$ and
$v_1=y\chi_{E_{\alpha_0}^{\gamma_0}}$, so that the functions $u_1$  and
$a_1$ are  supported  by the   same set. In this case we accordingly
choose $\alpha_1=\gamma_0$ and $\beta_1=\beta_0$.)

In accordance with our construction the support sets        
of  $u_1$ and $v_1$ are the disjoint  sets $E_{\alpha_0}^{\gamma_0}$ and
$E_{\gamma_0}^{\beta_0}$ respectively, which are left invariant by $A$.
The  same   disjoint  sets are the support sets of the elements   
$a_1$  and $b_1$.  This implies (in view of the equality $x=a_1\oplus b_1$) 
that 
$y=Ax =Aa_1\oplus Ab_1$, and therefore
$$
Aa_1=u_1\quad{and}\quad Ab_1=v_1\,.
$$

In the next lemma, we present some simple estimates on the norms of
$a_1$ and $b_1$.

\begin{lemma}\label{l:norm.est}  %\marginnote{l:norm.est}
For the functions $a_1$ and $b_1$ introduced above, we have:
\[
\textstyle
c\frac{\|y\|}{\|A\|}\le \|a_1\| \le c \|x\|\qquad\hbox{and}\qquad
 \|b_1\|  \le c_1 \|x\|\,, 
\]
where 
$\displaystyle c=1/{\root p \of 2}\,$ 
and 
$\displaystyle  c_1=\bigl[1-\bigl( \frac{c \|y\|}{\|A\|\cdot
\|x\|}\bigr)^p\bigr]^{1/p}>0$.  
\end{lemma}
\begin{proof}
Since $a_1\oplus b_1=x$ and $\|a_1\| \le\|b_1\|$, the $p$-additivity of         
the norm yields 
$$
2\|a_1\|^p\le \|a_1\|^p +\|b_1\|^p  =\|x\|^p\,,  
$$
whence  $\|a_1\|\le c\|x\|$.

From  $u_1=Aa_1$ we have  $\|u_1\|\le \|A\|\|a_1\|$.
So, taking into account that (in view of Lemma~\ref{l:norm})
$\|u_1\|=\|u_2\|=c\|y\|$, we see that 
$$
\displaystyle
c\frac{\|y\|}{\|A\|}=\frac{\|u_1\|}{\|A\|}\le\|a_1\|\,.
$$

For the last inequality, note that
\begin{eqnarray*}
\|b_1\|^p  &=&\|x\|^p -\|a_1\|^p \\
&\le& 
\|x\|^p - \frac{c^p\|y\|^p}{\|A\|^p} \\
&=&
\|x\|^p \bigl[ 1- \frac{c^p\|y\|^p}{\|A\|^p \|x\|^p}\bigr] =c_1^p\|x\|^p\,,
\end{eqnarray*}
and the proof of the lemma is finished.
\end{proof}

The rest of the construction must be done inductively. For instance, 
at the next step we will apply the above described procedure to the functions
$u_1, a_1$ satisfying $u_1=Aa_1$ and to the interval $[\alpha_1,\beta_1]$.  That
is, we take
$u_1$ for $y$ and 
$a_1$ for $x$ and we repeat the same procedure, keeping in mind that the    
support set of either of these two functions lies in
$E_{\alpha_1}^{\beta_1}$. 

Afterwards, we will have $u_1=u_2\oplus v_2$ with $\|u_2\| =\|v_2\|
=c\|u_1\|$ and with the support sets of these new functions
also invariant under $A$. 
Next we will have $a_1=a_2\oplus b_2$ with $\|a_2\|\le \|b_2\|$
and    
$Aa_2=u_2$, $Ab_2=v_2$ and with the corresponding estimates on the           
norms of $a_2,b_2$. The precise details of this inductive construction     
can be formulated as follows.

Assume that we have already constructed the functions $u_k$,$v_k$,$a_k$ and
$b_k$ and scalars  $\alpha_{k-1}<\gamma_{k-1} <\beta_{k-1}$   satisfying the   
following conditions: 
\begin{eqnarray*}
\|u_k\|&=&\|v_k\|=c\|u_{k-1}\|=c^k\|y\|\\
 u_{k-1}&=& u_k\oplus v_k \\
{\rm Supp}(u_k) &\subseteq& E_{\alpha_{k-1}}^{\gamma_{k-1}}\\
{\rm Supp}(v_k) & \subseteq& E_{\gamma_{k-1}}^{\beta_{k-1}}\\
 a_k &=& a_{k-1}\chi_{E_{\alpha_{k-1}}^{\gamma_{k-1}}} \\
 b_k &=& a_{k-1}\chi_{E_{\gamma_{k-1}}^{\beta_{k-1}} } \\
\hskip 4cm c^k{\frac{\|y\|}{\|A\|}}&\le&\|a_k\| \le c^k\|x\| \hskip 5.7cm (1) \\
\hskip 4cm \|a_k\|&\le&\|b_k\|\le c_1 c^{k-1}\|x\|   \hskip 5cm (2)
\end{eqnarray*}

For this choice of $a_k$ and $b_k$ we let  $\alpha_k=\alpha_{k-1}$             
and $\beta_k=\gamma_{k-1}$. 

Now we are ready to describe the induction step            
to produce $u_{k+1}$, $v_{k+1}$, $a_{k+1}$ and $b_{k+1}$, and the scalars
$\alpha_{k+1}$ and $\beta_{k+1}$.
Namely, to the elements $u_k$, $a_k$, satisfying  $u_k =Aa_k$, we apply the
very first step described in detail above. As a consequence, we  find first
the   scalar $\gamma_k \in (\alpha_k,\beta_k)$ such that
the functions $u_k\chi_{ E_{\alpha_k}^{\gamma_k}}$
and $u_k\chi_{ E_{\gamma_k}^{\beta_k}}$
have the same norm
$$
\|u_k\chi_{ E_{\alpha_k}^{\gamma_k}}\|
      =\|u_k\chi_{ E_{\gamma_k}^{\beta_k}}\|                                    
      =c\|u_k\|.  
$$
Next we consider the functions $a_k\chi_{ E_{\alpha_k}^{\gamma_k}}$         
and $a_k \chi_{ E_{\gamma_k}^{\beta_k}}$ and denote by $a_{k+1}$     
the one with the smaller norm---if both have the same norm, $a_{k+1}$          
can be either one. The other function is denoted  by $b_{k+1}$. Without       
loss of generality, we can  assume that
$a_{k+1}= a_k\chi_{ E_{\alpha_k}^{\gamma_k}}$. Subsequently, we let
$\alpha_{k+1} = \alpha_k$ and $\beta_{k+1} =\gamma_k$.
(Recall however, that if  
$\|a_k\chi_{ E_{\gamma_k}^{\beta_k}}\| < 
\|a_k\chi_{E_{\alpha_k}^{\gamma_k}}\|$, then 
$a_{k+1}= a_k\chi_{ E_{\gamma_k}^{\beta_k}}$, and accordingly
$\alpha_{k+1} = \gamma_k$ and $\beta_{k+1} =\beta_k$.)

We are ready to verify now that the functions $a_{k+1}$ and $b_{k+1}$
satisfy the desired estimates.

\begin{lemma}\label{l:induct}  %\marginnote{l:induct}

The functions  $a_{k+1},b_{k+1}$ constructed above satisfy the following
inequalities:
$$
 c^{k+1}\frac{\|y\|}{\|A\|} \le \|a_{k+1}\| \le
c^{k+1}\|x\|\qquad\hbox{and}\qquad
\|b_{k+1}\| \le c_1c^k\|x\|\,.
$$
\end{lemma}

\begin{proof}
By Lemma~\ref{l:norm.est} we have $\|a_{k+1}\|\le c\|a_k\|$.
This and the right 
inequality  in (1) imply that 
$\|a_{k+1}\|\le c^{k+1}\|x\|$. The equalities $Aa_{k+1}= u_{k+1}$  and
$\|u_k\|=c^k\|y\|$ imply 
$$
\|a_{k+1}\| \ge \frac{\|u_{k+1}\|}{\|A\|}=
c\frac{\|u_k\|}{\|A\|} =       
c^{k+1} \frac{\|y\|}{\|A\|}\,.
$$

Finally, we use the identity $a_{k+1}\oplus b_{k+1}= a_k$
and again the above estimate  $\|a_k\|\le c^k\|x\|$ 
to get:
\begin{eqnarray*}
\|b_{k+1}\|^p &=& \|a_k\|^p - \|a_{k+1}\|^p \le 
(c^{k} \|x\|)^p -\|a_{k+1}\|^p \\
&\le&
c^{kp} \|x\|^p - \bigl(c^{k+1} \frac{\|y\|}{\|A\|}\bigr)^p\\
&=&
c^{kp} \|x\|^p\bigl[1 - \bigl(\frac{c\|y\|}{\|A\|\|x\|}\bigr)^p\bigr]=
c_1^p c^{kp} \|x\|^p\,.
\end{eqnarray*}
This implies $\|b_{k+1}\| \le c_1 c^k \|x\|$, as desired.
\end{proof}

Using the sequence $\{b_k\}$ and the estimates obtained so far,
we can finally produce  a sequence $\{e_k\}$ satisfying the               
properties required  in the proof of Theorem~\ref{t:oc.const}.

\begin{lemma}\label{l:final}    %\marginnote{l:final}
If $e_n=\frac{b_n}{\|b_n\|}$, then the sequence $\{e_n\}$ satisfies the
following properties:
\begin{itemize}
\item[(i)]
$\|e_n\| =1$ for each $n$,
\item[(ii)]
$\{Ae_n\}$ is a disjoint sequence, and
\item[(iii)]
$\|Ae_n\| \ge \delta$ for each $n$ and for some $\delta >0$.
\end{itemize}
\end{lemma}
\begin{proof}
Since, by their definition, the vectors $b_n$ are pairwise disjoint
and have the  sets  $E_{\gamma_n}^{\beta_{n-1}}$ (which are disjoint           
and invariant under our operator $A$) as their support sets, we see           
that the vectors
$Ab_n$, $n=1,2,\ldots$, are also pairwise disjoint. Now
recalling that $Ab_n=v_n$ and using the right 
inequality in (2)       
we can easily estimate $\|Ae_n\|$:  
\begin{eqnarray*}
\textstyle
\|Ae_n\|&=&\frac{\|Ab_n\|}{\|b_n\|} =\frac{\|v_n\|}{\|b_n\|} 
 =c^n\frac{\|y\|}{\|b_n\|}\\
\textstyle
&\ge&
\frac{c^n}{c_1c^{n-1}} \frac{\|y\|}{ \|x\|} =
\frac{c}{c_1}\cdot\frac{\|y\|}{\|x\|}\,. 
\end{eqnarray*}
This completes the proof.      
\end{proof}

\end{document}